\documentclass[12pt]{amsart}
\usepackage{amsmath,amssymb,amscd,amsfonts}
\newtheorem{prop}{Proposition}[section]
\newtheorem{thm}[prop]{Theorem}
\newtheorem{cor}[prop]{Corollary}

\newtheorem{lem}[prop]{Lemma}
\def\demo#1{\medskip\noindent{\bf#1\enspace}}

\def \calT {{\mathcal T}}
\def \R {{\mathbb R}}
\def \S {{\mathcal I}}
\def\QED{{\hfill $\square$ \medskip}}
        
        \def\zed{{\mathbb Z}}
        \def\Z{{\zed}}
        \def\z{{\bf z}}
        \def\v{{\bf v}}
        \def\k{{\bf k}}
        \def\t{{\bf t}}
        \def\i{{\bf i}}
        \def\L{{\bf L}}
        
        \def\r{{\bf r}}
        \def\red{{\mathbb R}}
        
        \def\qed{{\mathbb Q}}
        
\input epsf.tex
\begin{document}

\title{When shape matters: deformations of tiling spaces}

\author{Alex Clark and Lorenzo Sadun}

\address{Alex Clark: Department of Mathematics, University of North Texas,
             Denton, Texas 76203}
\email{AlexC@unt.edu}

\address{Lorenzo Sadun: Department of Mathematics, The University of Texas at Austin,
 Austin, TX 78712-1082 U.S.A.}
\email{sadun@math.utexas.edu}   



\begin{abstract}
We investigate the dynamics of tiling dynamical systems and their
deformations.  If two tiling systems have identical combinatorics,
then the tiling spaces are homeomorphic, but their dynamical
properties may differ.  There is a natural map $\S$ from the parameter
space of possible shapes of tiles to $H^1$ of a model tiling space,
with values in $\R^d$. Two tiling spaces that have the same image
under $\S$ are mutually locally derivable (MLD).  When the
difference of the images is ``asymptotically negligible'', then the
tiling dynamics are topologically conjugate, but generally not MLD. For
substitution tilings, we give a simple test for a cohomology class to
be asymptotically negligible, and show that infinitesimal deformations
of shape result in topologically conjugate dynamics only when the
change in the image of $\S$ is asymptotically negligible.  Finally, we
give criteria for a (deformed) substitution tiling space to be
topologically weakly mixing.
\end{abstract}

\maketitle

\markboth{Clark-Sadun}{When shape matters: Deformations of tiling
spaces}
\newpage

\section{Introduction and Statement of Results}

A tiling is described by a combination of combinatorial data (which
tiles meet which others) and by geometric data (the shape and location
of each tile).  Tilings with the same combinatorics may have different
geometry.   Compare figure 1, which shows a patch of
a Penrose tiling, to figure 2, which shows a corresponding
patch of a combinatorially identical but geometrically different 
tiling.  The translation group acts on the continuous hulls of both tilings.
How do the dynamics compare? 

\begin{figure}
\vbox{\hskip 1 in \epsfxsize=3truein\epsfbox{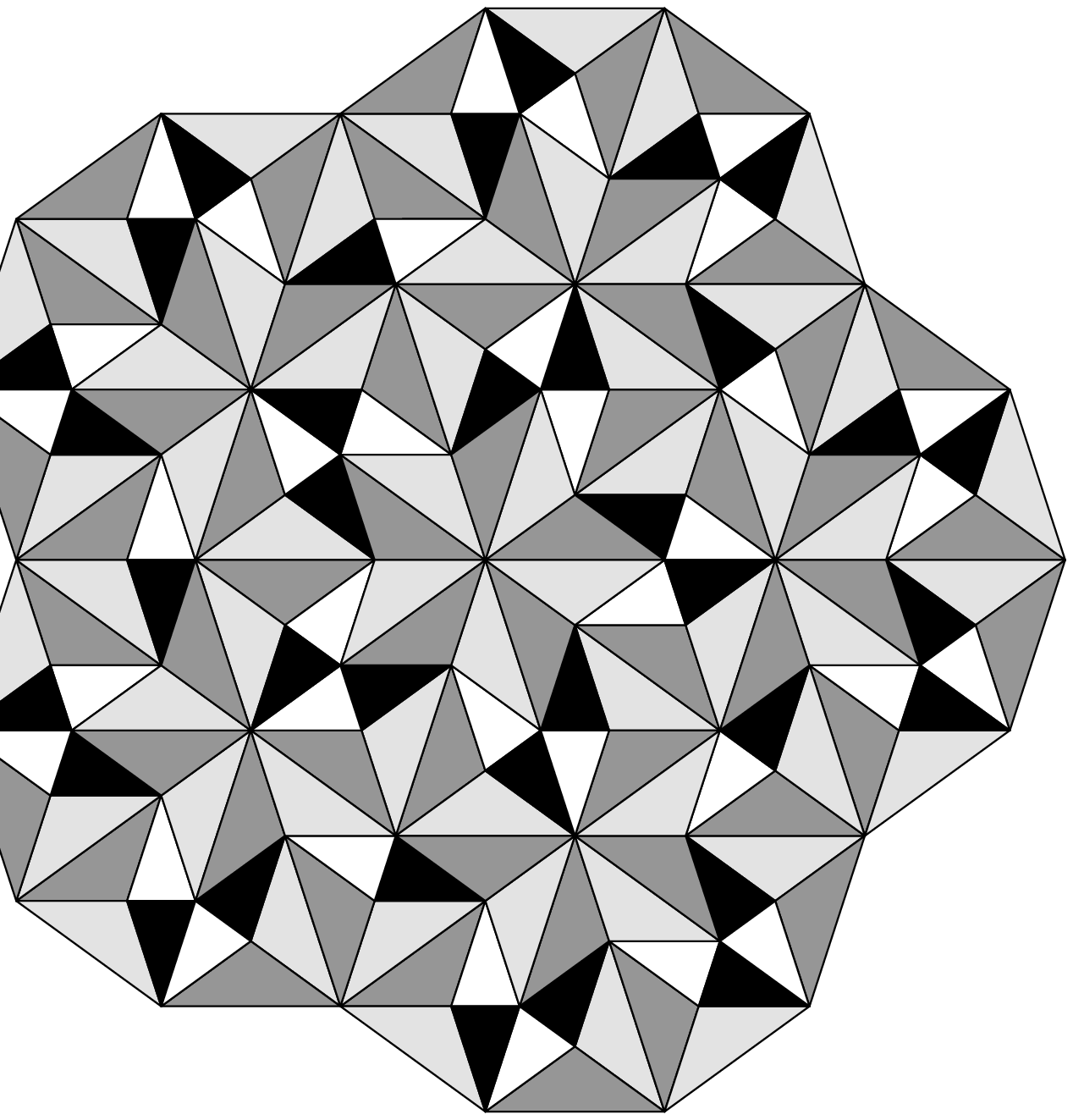}}
\caption{A patch of a Penrose tiling}
\end{figure}

\begin{figure}
\vbox{\hskip 1 in \epsfxsize=3truein\epsfbox{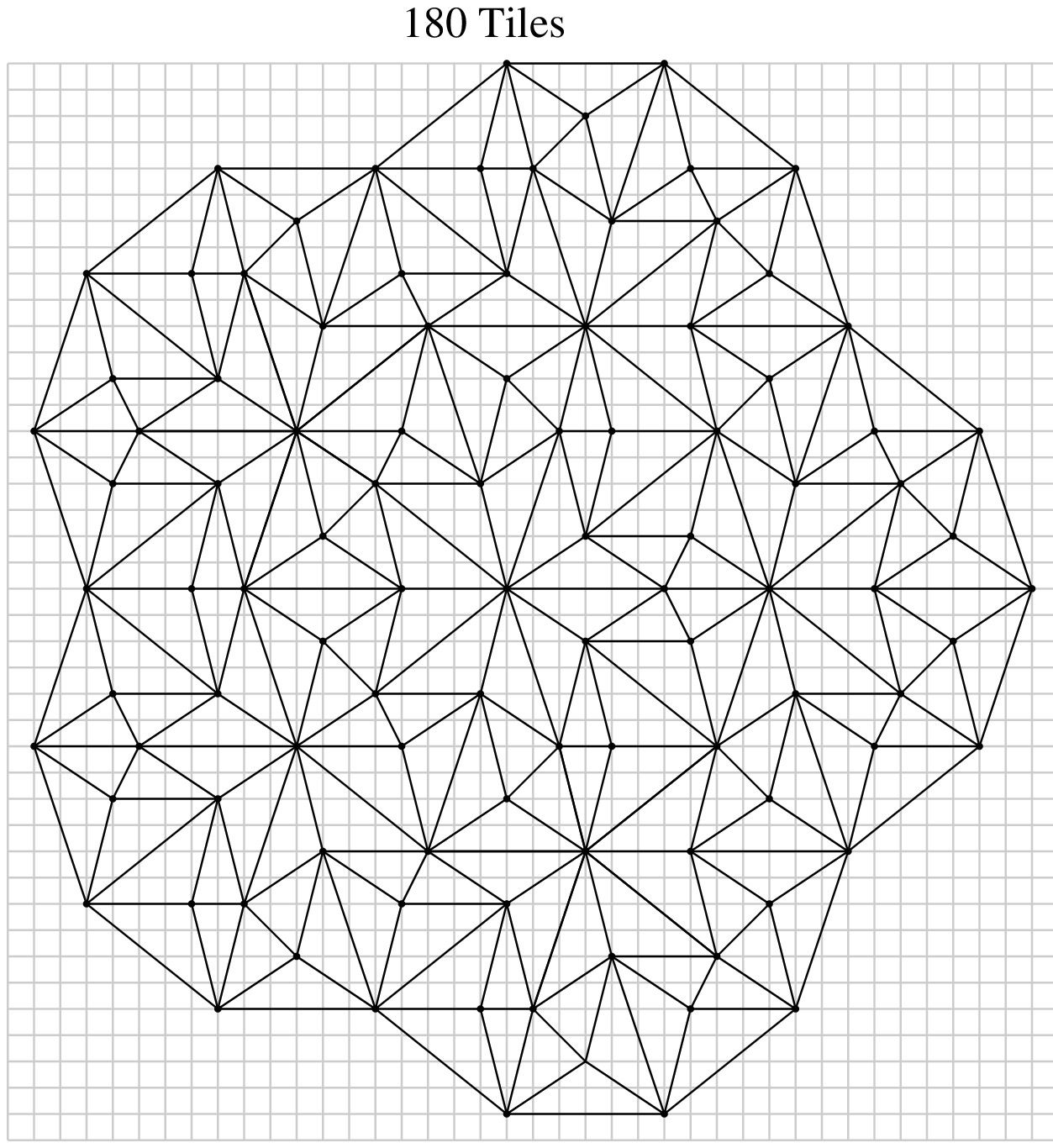}}
\caption{A corresponding patch of a deformed Penrose tiling}
\end{figure}

To avoid trivialities, we shall only consider \emph{nonperiodic}
tilings in this paper. Moreover, we shall assume that the tiles are
polyhedra that meet full-face to full-face, and that there are only a
finite set of tile types (a.k.a.~prototiles) up to translation.
Although these conditions may seem to be restrictive, any tiling that
has finite local complexity with respect to translations is mutually
locally derivable (MLD, see definition below) to a tiling meeting
these conditions, via derived Voronoi tilings \cite{Pr}. We also 
assume that our tiling spaces are minimal.  This is equivalent to the 
condition that each patch of each tiling appear in all other tilings
with bounded gaps. 

Deformations of tiling spaces were considered in \cite{SW}.  As
detailed there, to describe the shape and size of each tile, one must
specify the displacement vector that corresponds to each edge of each
prototile. If, somewhere in the tiling, translations of two prototiles
meet along a common edge, then those two edges must have the same
associated displacement vector.  Moreover, the sum of the displacement
vectors around each tile must be zero.  Each deformation thus
corresponds to a solution of a homogeneous linear system of equations
subject to additional conditions which are open in the solution space:
the sequence of edges around a tile does not cross itself and has
winding number +1 and bounds a $2$-cell of a tile.  If $d>2$ we also
require the edges of a $2$-cell to be coplanar, and in general for all
the edges of a $k$-cell of a tile to lie in a $k$-dimensional
hyperplane of $\R^d$.

Put another way, the shapes of the prototiles are determined by a
function $f$ from the set of prototile edges (modulo certain
identifications) to $\R^d$.  This function may be viewed as a
1-cochain on the CW-complex obtained by taking the disjoint union of
all the prototiles, modulo the identification of edges (and other
structures of dimension less than $d$) where tiles can meet.  This is
precisely the complex $\Gamma$ considered by Anderson and Putnam in
their computation of topological invariants of substitution tiling
spaces \cite{AP}.  (The Anderson-Putnam complex is defined for all
translationally finite tilings, although it has been applied primarily
to substitution tilings. See \cite{G} for applications in a more
general context.)  That is, $f \in C^1 (\Gamma, \R^d) =Hom(C_1(\Gamma), \R^d)$,
and we may think of ``shape space'' as a subset of
$C^1(\Gamma, \R^d)$. (See, e.g., \cite{M} for background on (co)homology.)

We denote vectors by bold-face roman letters, 
tiling spaces by $\calT$ with appropriate subscripts, 
individual tilings by roman letters from the end of the alphabet (usually
$x$ or $y$), and shape parameters by $f$ or $g$. The translate of the
tiling $x$ by $\z \in \R^d$ is denoted $x-\z$.
   
Fix a tiling space $\calT$, and consider deformations of this space.
For any choice $f$ of edge vectors meeting the above requirements, we
can consider tilings whose tiles are described by these edge vectors,
but whose combinatorics (which tiles meet and how) are the same as the
original tilings in $\calT$.  Let $\mathcal{T}_{f}$ denote the space
of such deformed tilings, on which $\R^d$ acts naturally by
translation.  The primary focus of this paper is the extent to which
the dynamical system $\mathcal{T}_{f}$ depends on the function $f$. If
all we care about is the topological space, a theorem of \cite{SW}
shows that it doesn't:

\begin{thm}\cite{SW}
\label{homeo}Let $f,g \in C^1(\Gamma, \R^d)$ be admissible shape functions
for a fixed tiling space.
Then $\calT_{f}$ and $\calT_g$
are homeomorphic.
\end{thm}

But what about dynamics? In a previous paper \cite{CS} we showed how
1-dimensional substitution tiling space dynamics depend on the lengths
of the tiles.  Here we extend that analysis to higher dimensions, and
also recast our 1-dimensional results in terms of topological
invariants.
One natural question is when two tiling spaces have topologically conjugate
dynamics.  Another is when they have dynamics that are intertwined by a 
local map.  This involves the notion of {\em mutual local derivability},
first introduced in \cite{BSJ}. 

Let $x$ and $y$ be two tilings, possibly with different sets of
prototiles. The tiling $y$ is said to be {\em locally derivable} from $x$
if there exists a length $R$ such that, if $\z_1, \z_2 \in \R^d$ and
$x-\z_1$ agrees with $x-\z_2$ on a ball of radius $R$ around the
origin, then $y-\z_1$ agrees with $y-\z_2$ on a ball of radius 1
around the origin.  In other words, the type and exact placement of
the tile at a point $\z$ in $y$ depends only on the patch of radius
$R$ around $\z$ in $x$.  If $y$ is locally derivable from $x$ and $x$ is
locally derivable from $y$, then $x$ and $y$ are said to be mutually locally
derivable (MLD) tilings.

The MLD concept extends to tiling spaces.  If $x$ and $y$ are MLD
tilings, then the closure of the translational orbit of $x$ is
topologically conjugate to the closure of the translation orbit of
$y$, via a conjugacy that takes $x$ to $y$, and thus depends only on
local data.  More generally, we say that two tiling spaces $\calT$ and
$\calT'$ are MLD if there exists a topological conjugacy $\phi: \calT
\mapsto \calT'$ such that, for some $R$, the patch of size 1 around
the origin in $\phi(x)$ can be determined exactly from the patch of
size $R$ around the origin in $x$, and the patch of size 1 in $x$ can
be determined exactly from the patch of size $R$ in $\phi(x)$.  This
is a natural generalization of the concept of ``sliding block codes''
for subshifts.  However, while all continuous maps of subshifts are
sliding block codes (see, e.g., \cite{LM}), there exist topologically
conjugate tiling spaces that are not MLD \cite{RS,Pe,CS}.

In Section 2, we construct a natural map $\S$ from the space of shape
parameters to the Cech cohomology group $H^1(\calT,\R^d)$ 
($=H^1(\calT, \Z) \otimes \R^d$) and show
that the dynamical properties of $\calT_f$ depend only on $\S(f)$:

\begin{thm}[Theorem \ref{mld2}]\label{mld}
Suppose $f$ and $g$ are admissible shape parameters with $\S(f)=\S(g)$.
Then $\calT_f$ and $\calT_g$ are mutually locally derivable.
\end{thm}

Furthermore, we define a condition on $H^1(\calT, \R)$ (and by extension,
$H^1(\calT,\R^d)$) called ``asymptotic negligibility''.  Roughly
speaking, a deformation is asymptotically negligible if, up to an
arbitrarily small error, it does not change the return vectors of
large patches in a tiling.  We then prove

\begin{thm}[Theorem \ref{negligible2}]\label{negligible}
Let $f,g \in C^1(\Gamma, \R^d)$ be admissible shape functions
for a fixed tiling space. If $\S(f)-\S(g)$ is asymptotically negligible, 
then there is a topological conjugacy between 
$\calT_{f}$ and $\calT_g$.
\end{thm}

Rescaling a tiling, or more generally applying a linear transformation, 
does not change the qualitative dynamical properties of the tiling space,
such as mixing, minimality, or diffractivity. We therefore consider when
a tiling space $\calT_f$ is topologically conjugate to a linear transformation
of a tiling space $\calT_g$. Since linear transformations act naturally on
the $\R^d$ factor of $H^1(\calT, \R^d)$ we have:

\begin{cor}\label{neg-cor}
The tiling spaces $\calT_f$ and $\calT_g$
are conjugate up to linear transformation if, for some linear transformation
$L$, $\S(g) - L (\S(f))$ is asymptotically negligible. 
\end{cor}

In Section 3 we specialize to nonperiodic tilings made from a primitive 
substitution, or a primitive substitution-with-amalgamation. The resulting
spaces are orbit closures of self-similar or pseudo-self-similar tilings, 
respectively, in the terminology of \cite{So1,So2,PS}.   
An essential object of study is the action of the substitution on the
tiling space itself, and therefore on $H^1(\calT,\R^d)$. We decompose
$H^1(\calT, \R^d)$ into (generalized) eigenspaces of the substitution
operator.  Let $S(\calT)$ be the span of the (generalized) eigenvectors with
eigenvalue strictly less than one in magnitude. We show that the set of
asymptotically negligible classes is precisely $S(\calT)$, and prove the 
following local converse to Theorem \ref{negligible}:

\begin{thm}[Theorem \ref{necessary0}]\label{converse}
For each shape parameter $f$ there is a neighborhood $U_f$ of $\S(f)$
such that, if $\S(g) \in U_f$ and $\calT_f$ and $\calT_g$ are 
topologically conjugate,
then $\S(g)-\S(f) \in S(\calT)$. 
\end{thm}

If there is a pseudo-self-similar tiling in $\calT_f$,
then $\S(f)$ is a Perron-Frobenius eigenvector of the substitution
applied to $H^1(\calT, \R^d)$. The Perron-Frobenius eigenspace, denoted
$PF(\calT)$, is $d^2$ dimensional, and precisely equals $L \S(f)$, where
$L$ ranges over all linear transformations of $\R^d$.  Corollary \ref{neg-cor} 
has the following simple restatement for substitution tilings:

\begin{cor} \label{neg-cor-sub}
If $\calT_f$ contains a pseudo-self-similar tiling, and if $\S(g)
\in PF(\calT) \oplus S(\calT)$, then $\calT_g$ is topologically conjugate
to a linear transformation applied to $\calT_f$. 
\end{cor}

In particular, if $H^1(\calT,\R^d)= PF(\calT) \oplus S(\calT)$, then all
choices of shape give rise to conjugate dynamics, up to linear transformation.
This case is the natural generalization, to higher dimensions, of Pisot
substitutions. 



For the space of Penrose tilings, $H^1(\calT,\Z) = \Z^5$ \cite{AP},
so $H^1(\calT,\R^2)$ is 10-dimensional.  
The eigenvalues of the substitution are the golden mean $\tau$,
with multiplicity 4, $1-\tau$ with multiplicity 4, and $-1$ with
multiplicity 2. However, the $-1$ eigenspace is odd under rotation by
180 degrees, while the $\tau$ and $1-\tau$ eigenspaces are even \cite{ORS}.
Any choice of shapes that preserves the (statistical) 2-fold rotational
symmetry of the Penrose tiling must not involve the $-1$ eigenspace.
As a result, the tiling space constructed from the tiling of figure 2
is topologically conjugate to a linear transformation applied to the
undistorted Penrose tiling space.

The techniques of Section 3 require recognizability, and hence draw heavily 
on Solomyak's generalization \cite{So2} of the work of Moss\'e 
on 1-dimensional subshifts \cite{M1, M2}. 

In Section 4 we study the spectra of substitution tilings and their
deformations. We provide general criteria for the existence of point
spectrum of translations acting on $\calT_f$, similar in spirit to
criteria found in \cite{So1} for self-similar tilings and to criteria
in \cite{CS} for 1-dimensional substitution tiling spaces and their
deformations.  We also provide constraints on the form of that
spectrum, in terms of the aforementioned decomposition of $\S(f)$ into
eigenspaces of the substitution operator.  In particular, we show:

\begin{thm}[Theorem \ref{somesmall1}]\label{mixing}
If $H^1(\calT,\R^d) \neq PF(\calT) \oplus S(\calT)$, then for a
generic choice of shape parameter $f$, $\calT_f$ is topologically
weakly mixing.
\end{thm}

Finally, in Section 5 we revisit the problem in one dimension, and recast
the results of \cite{CS} in topological terms. 

\section{The Map $\S$}

To define the map $\S$ we first recall the inverse limit structure of 
tiling spaces, developed by Anderson and Putnam \cite{AP} 
for substitution tilings,
and generalized by G\"{a}hler \cite{G} to apply to all translationally finite
tilings. (See \cite{BBG} for an alternate approach, and \cite{ORS, Sa, BG}
for further generalizations). 

We first construct a complex $\Gamma$ by taking the
disjoint union of prototiles in the tiling, modulo identification of edges
where tiles can meet.  If somewhere in a tiling, edge $i$ of a tile of type 
$A$ is coincident with edge $j$ of a tile of type $B$, then we identify
edge $i$ of $A$ with edge $j$ of $B$ in the complex.  (If $d>2$ we also
identify coincident faces and other structures of dimension up to $d-1$.)

We may also rewrite the tiling in terms of collared tiles, labeling each
tile $t$ by the patch consisting of all tiles that touch $t$. Applying
the construction of the previous paragraph to the collared tiles gives a 
complex $\Gamma^{(1)}$.  Collaring the collared tiles and applying
the construction gives a complex $\Gamma^{(2)}$, and more generally applying
the construction to $k$-times collared tiles gives a complex $\Gamma^{(k)}$. 

There is a natural ``forgetful'' map $\alpha_0$ from $\Gamma^{(1)}$ to
$\Gamma$ that simply ignores the collaring, and likewise a forgetful
map $\alpha_{k}: \Gamma^{(k+1)} \mapsto \Gamma^{(k)}$.  The inverse limit
of the sequence of maps and spaces
\begin{equation}
\Gamma = \Gamma^{(0)} \overset{\alpha_0}{\leftarrow} \Gamma^{(1)} 
\overset{\alpha_1}{\leftarrow} \Gamma^{(2)} \cdots
\end{equation}
is isomorphic to the tiling space $\calT$ \cite{G}. 
A point in $\Gamma$ tells how to
place a tile around the origin.  A point in $\Gamma^{(1)}$ tells how
to place a collared tile, i.e., a tile and its nearest neighbors. As
$k$ increases, the points in $\Gamma^{(k)}$ tell how to place larger and
larger patches around the origin, and the entire sequence $(x_0, x_1, \ldots)$
with $x_i \in \Gamma^{(i)}$ and $x_i = \alpha_i(x_{i+1})$ tells how to 
place a complete tiling. 

We have already seen that the a shape parameter $f$ is a vector-valued
1-cochain in $\Gamma$. In fact, for any tile, $\delta f(t) =
f(\partial t)=0$, since this is the sum of the edge vectors around the
tile $t$. Thus $f$ is a cocycle, and defines a cohomology class $[f]
\in H^1(\Gamma, \R^d)$.  Let $r_k: \calT \mapsto \Gamma^{(k)}$ be the
projection of the tiling space to the $k$th approximant.  Define
$\S(f) = r_0^* [f] \in H^1(\calT,\R^d)$.

\begin{thm}[Theorem \ref{mld}]
Suppose $f$ and $g$ are admissible shape parameters with $\S(f)=\S(g)$.
Then $\calT_f$ and $\calT_g$ are mutually locally derivable.
\label{mld2}
\end{thm}

\demo{Proof.} We show how to construct a tiling in $\calT_g$ from the 
corresponding tiling in $\calT_f$ by a process that is completely local. 
The reverse process is of course similar. 

Let $\pi_k = \alpha_0 \circ \alpha_1 \circ \cdots \circ \alpha_{k-1}:
\Gamma^{(k)} \mapsto \Gamma$. If $\S(f)=\S(g)$, then for some finite $k$,
$\pi_k^*[f] - \pi_k^*[g] = 0$, so $\pi_k^*f - \pi_k^*g = \delta \beta$
for some $\beta \in C^0(\Gamma^{(k)}, \R^d)$. Now each vertex $\v$ in a
tiling in $\calT_f$ maps to a unique vertex in $\Gamma^{(k)}$, the map
being determined by a ball of size $(k+1) A$ around $\v$, where $A$ is the
diameter of the largest prototile.  Moving each vertex $\v$ by $-\beta(\v)$,
and linearly interpolating the edges between vertices, converts a tiling
in $\calT_f$ to a tiling in $\calT_g$. \QED

For a tiling $x \in \calT$, a \textit{recurrence} is an ordered pair
$(\z_1,\z_2)$ of points in $x$ such that $\z_1$ is a point in a
tile and $\z_2$ is the corresponding point of a translate of that
tile and such that there exist balls around $\z_1$ and $\z_2$ that
agree (up to translation by $\z_2-\z_1$, of course).  If $r$ is the
supremum of the radii of the balls around $\z_1$ and $\z_2$ that
agree, then we say the recurrence has \textit{size} $r$.  Each path
along edges from $\z_1$ to $\z_2$ projects to a closed loop in $\Gamma$,
and hence to a closed chain in $C_1(\Gamma)$.  Different paths from
$\z_1$ to $\z_2$ correspond to homologous chains. The class in
$H_1(\Gamma)$ of a recurrence is called a \textit{recurrence class}.
Note that recurrences of size greater than $(k+1)A$, where $A$ is the
diamater of the largest prototile, also project to closed paths in
$\Gamma^{(k)}$ and define classes in $H_1(\Gamma^{(k)})$.  Since the
tiling space $\calT$ is assumed to be minimal, the set of recurrence
classes is the same for every tiling in the space, and so we can speak
of the recurrence classes of the tiling \textit{space}.

An element $\eta$ of $H^1(\Gamma^{(k)}, \R)$, or of $H^1(\Gamma^{(k)},
\R^d)$, is said to be \textit{asymptotically negligible} if, for each
$\epsilon>0$ there exists a constant $R_\epsilon$ such that $\eta$,
applied to any homology class that can be represented by a 
recurrence of size greater than $R_\epsilon$, is less
than $\epsilon$ in magnitude.  A class in $H^1(\calT, \R^d)$ is
asymptotically negligible if it is the pullback of an asymptotically
negligible class in $H^1(\Gamma^{(k)}, \R^d)$ for some finite $k$. It
is clear that any linear combination of asymptotically negligible
classes is asymptotically negligible, so these classes form a subspace
of $H^1(\calT, \R^d)$, denoted $N(\calT)$.

\begin{thm}[Theorem \ref{negligible}]\label{negligible2}
Let $f,g \in C^1(\Gamma, \R^d)$ be admissible shape functions
for a fixed tiling space. If $\S(f)-\S(g)$ is asymptotically negligible, 
then there is a topological conjugacy $\phi: \calT_{f} \mapsto \calT_g$.
\end{thm}

\demo{Proof.} We construct the conjugacy $\phi$ in stages. First pick
a reference tiling $x \in \calT_f$ that has a vertex at the origin. 

For every vertex $\v$ in $x$ there is a path $p_\v$ from the origin to $\v$
along edges, and the location of $\v$ is precisely $f(p_\v)$.  In $\phi(x)$
we place the corresponding vertex at $g(p_\v)$. The path $p_\v$ is not uniquely
defined, but different choices differ by boundaries, so the values of
$f(p_\v)$ and $g(p_\v)$ are uniquely determined.  Once the location of
the vertices of $\phi(x)$ are specified, constructing the edges and tiles is
straightforward.  

This defines $\phi(x)$. For $\z \in \R^d$, let $\phi(x-\z)=\phi(x)-\z$.
It remains to show that $\phi$ is uniformly continuous on the orbit of $x$, and
hence can be extended to all of $\calT$. 

Suppose that $\z_1$ and $\z_2$ are vertices in $x$ such that $x-\z_1$
and $x-\z_2$ agree on a large ball around the origin.  Note that $\phi(x) -
g(p_{\z_1})$ agrees exactly with $\phi(x) - g(p_{\z_2})$ on a large
ball around the origin.  Since $f-g$ is asymptotically negligible,
$f(p_{\z_2})-f(p_{\z_1})$ is very close to $g(p_{\z_2})-g(p_{\z_1})$,
so $\phi(x-\z_1)=\phi(x) - f(p_{\z_1})$ agrees with $\phi(x-\z_2) =
\phi(x) - f(p_{\z_2})$ on a large ball around the origin, up to a
small translation, of size $f(p_{\z_2})-f(p_{\z_1}) - \left (
g(p_{\z_2})-g(p_{\z_1}) \right )$.  Since this translation can be made
arbitrarily small by making the recurrence of sufficiently large size
(independent of the choice of $\z_1$ and $\z_2$), $\phi$ is uniformly
continuous on the orbit of $x$.

To see that $\phi$ is invertible, construct a semi-conjugacy $\phi': \calT_g
\mapsto \calT_f$ by the same procedure, with the roles of $f$ and $g$ reversed,
and with $\phi(x)$ as the reference tiling. Since $\phi'(\phi(x))=x$, it
is clear that $\phi$ and $\phi'$ are inverses. \QED

Note that $\calT_f$ and $\calT_g$ are typically not MLD. If two
tilings in $\calT_f$ agree on a large ball around the origin, the
corresponding tilings in $\calT_g$ agree on a large ball \textit{up to a small
translation}. Only if $f-g$ vanishes on large recurrence classes do
they agree without translation. In that case, however, the following
theorem shows that $\S(f) = \S(g)$, and we are back in the situation
of Theorem \ref{mld}.

\begin{thm}\label{span-H_1}
Let $\beta \in H^1(\Gamma, \R)$.  If $\beta$ vanishes on all recurrence
classes of size greater than a fixed value $R$, then the pullback of $\beta$
to $H^1(\calT, \R)$ is zero. 
\end{thm}

\demo{Proof.}  We use Kellendonk and Putnam's P-equivariant cohomology 
\cite{K, KP},
which relates the real-valued cohomology of a tiling space to closed
and exact forms on a single tiling $x$, meeting some equivariance
conditions.  A differential form on the tiling $x$ (viewed as a
decorated copy of $\R^d$) is said to be P-equivariant if there is some
radius $r$ such that the value of the form at each point depends only
on the tiling in a ball of radius $r$ around that point.  That is, if
$\rho =\sum_{i_{1},\dots ,i_{k}}f_{i_{1},\dots ,i_{k}}dx_{i_{1}}\dots
dx_{i_{k}}$
 is an equivariant form with radius $r$, and if $x-\z_1$ and
$x-\z_2$ agree on a ball of radius $r$ around the origin, then
$f_{i_{1},\dots ,i_{k}}(\z_1)=f_{i_{1},\dots ,i_{k}}(\z_2)$. Kellendonk and Putnam prove a version of the de Rham
theorem: the closed P-equivariant $k$-forms on $x$, modulo $d$ of the
P-equivariant $(k-1)$-forms, is isomorphic to $H^k(\calT,\R)$.

Now consider our class $\beta$, and represent it as a closed P-equivariant
1-form $\rho$. Let $\gamma(\z) = \int_0^\z \rho$. It is clear that $\rho =
d\gamma$, and the condition that $\beta$ vanishes on large recurrence classes
implies that $\gamma$ is P-equivariant with radius $R$. Thus $\rho$ represents
the zero class in $H^1(\calT,\R)$.  \QED

\section{Substitution Tiling Spaces}
\label{sub}
We now specialize to substitution tiling spaces, where we can obtain 
sharper results.  In this section we identify the asymptotically negligible
cohomology classes and prove Theorem 1.5.

A substitution tiling system is determined by a substitution $\sigma$
from the set of prototiles to the set of finite patches, such that for
each prototile $t$, the tiles of $\sigma(t)$ do not overlap, and such that
their union is the rescaled prototile $\lambda t$, where
$\lambda$ is a fixed ``stretching factor''.  For example, in the
``chair'' tiling, there are four prototiles:  One is the L-shaped tile
of figure 3 and the others are this tile rotated by $\pi/2$, $\pi$ and 
$3 \pi/2$.
The substitution map is shown in figure 3.  We extend the map
$\sigma$ to patches, dilating the entire patch by a factor of $\lambda$
and replacing each dilated tile $\lambda t_i$ with $\sigma(t_i)$. 
For each prototile $t$ we have a sequence of
patches $t, \sigma(t), \sigma^2(t), \cdots$. A tiling $x$ of $\R^d$ by
prototiles is called admissible if for every finite patch $P$ of 
$x$, there exists a prototile $t$ and an integer $n$ such that $P$ is 
a translate of a subpatch of $\sigma^n(t)$. The substitution tiling space
$\calT$ is then the set of all admissible tilings, and $\sigma$ extends
to a continuous map $\calT \mapsto \calT$. 

\begin{figure}
\vbox{\hskip 1 in \epsfxsize=3truein\epsfbox{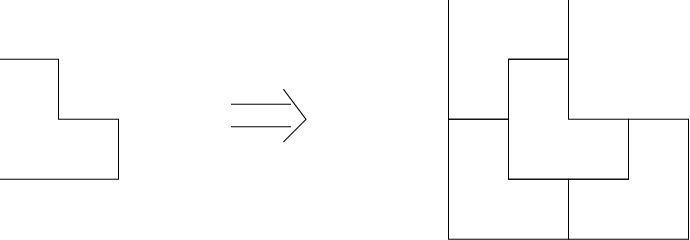}}
\caption{The chair tile and its substitution}
\end{figure}

We assume that the substitution is primitive.  That is, there exists
an integer $n$ such that, for any two (possibly identical) prototiles
$t_1$ and $t_2$, $\sigma^n(t_1)$ contains a copy of $t_2$. This
ensures that $\calT$ is minimal.  Furthermore, we assume that $\calT$
is non-periodic. There does not exist a tiling $x \in \calT$ and a
nonzero $\z \in \R^d$ with $x-\z = x$. By a theorem of Solomyak \cite{So2},
extending previous work by Moss\'e \cite{M1, M2}, this implies that $\sigma:
\calT \mapsto \calT$ is a homeomorphism. There then exists a ``recognition
length'' $D$ such that, if $x$ and $y$ are tilings that agree on a ball
of radius $r>D$ about a point $\z \in \R^d$, then $\sigma^{-1}(x)$ and
$\sigma^{-1}(y)$ agree on the tile(s) containing the point $\z/\lambda$.

In every substitution tiling space there exists a fixed point of some
power of the substitution.  Replacing $\sigma$ with a power of
$\sigma$ does not change the tiling space, so we can assume, without
loss of generality, that there is a tiling $x$ with $\sigma(x)=x$.
Such a tiling is called ``self-similar'', and has the property that
$\lambda x$ is a tiling by large tiles, each of which is a union of
tiles of $x$.  Some authors begin by defining self-similar
tilings by this property, and then take $\calT$ to be the closure of
the translational orbit of $x$.

Closely related to self-similar tilings are pseudo-self-similar
tilings.  A tiling $x$ is pseudo-self-similar if, for some scaling
factor $\lambda$, $\lambda x$ and $x$ are MLD. Natalie Priebe Frank
and Boris Solomyak \cite{PS} have shown that every pseudo-self-similar
planar tiling with polygonal tiles is MLD to a self-similar tiling.
However, the tiles of the self-similar tiling may not be polygonal;
rather, they may have fractile boundaries. Conversely, Natalie Priebe
Frank \cite{Pr} showed (in arbitrary dimension) that each self-similar
tiling (with tiles of arbitrary shape) is MLD to a pseudo-self-similar
tiling with polyhedral tiles meeting full-face to full-face.

Orbit closures of pseudo-self-similar tilings may also be viewed as
coming from a substitution-with-amalgamation.  A
substitution-with-amalgamation is the composition of a linear rescaling
and a local map.  The local map is a prescription for replacing
each rescaled prototile $\lambda t$ by a collection of tiles, such that 
the resulting tiles do not overlap and
do not leave any gaps. As with an ordinary substitution, $\sigma$
extends to a homeomorphism $\calT \mapsto \calT$. 
In a (strict) substitution tiling, each rescaled
prototile $\lambda t$ in $\lambda x$ is a union of tiles of $\sigma(x)$,
but in a substitution-with-amalgamation, the tiles of $\sigma(x)$ may stick in
and out of the tiles of $\lambda x$. 

An example of a substitution-with-amalgamation, for a
tiling by marked hexagonal tiles, is shown in figure 4.  Each small
hexagon is rescaled by a linear factor of two, and the resulting large 
hexagon is replaced by four small hexagons. The union of the four small 
hexagons is not the corresponding large hexagon, yet 
the small hexagons associated to {\it all} the large hexagons yield a tiling
of the plane, without gaps or overlaps. 

\begin{figure}
\vbox{\hskip 1 in \epsfxsize=2truein\epsfbox{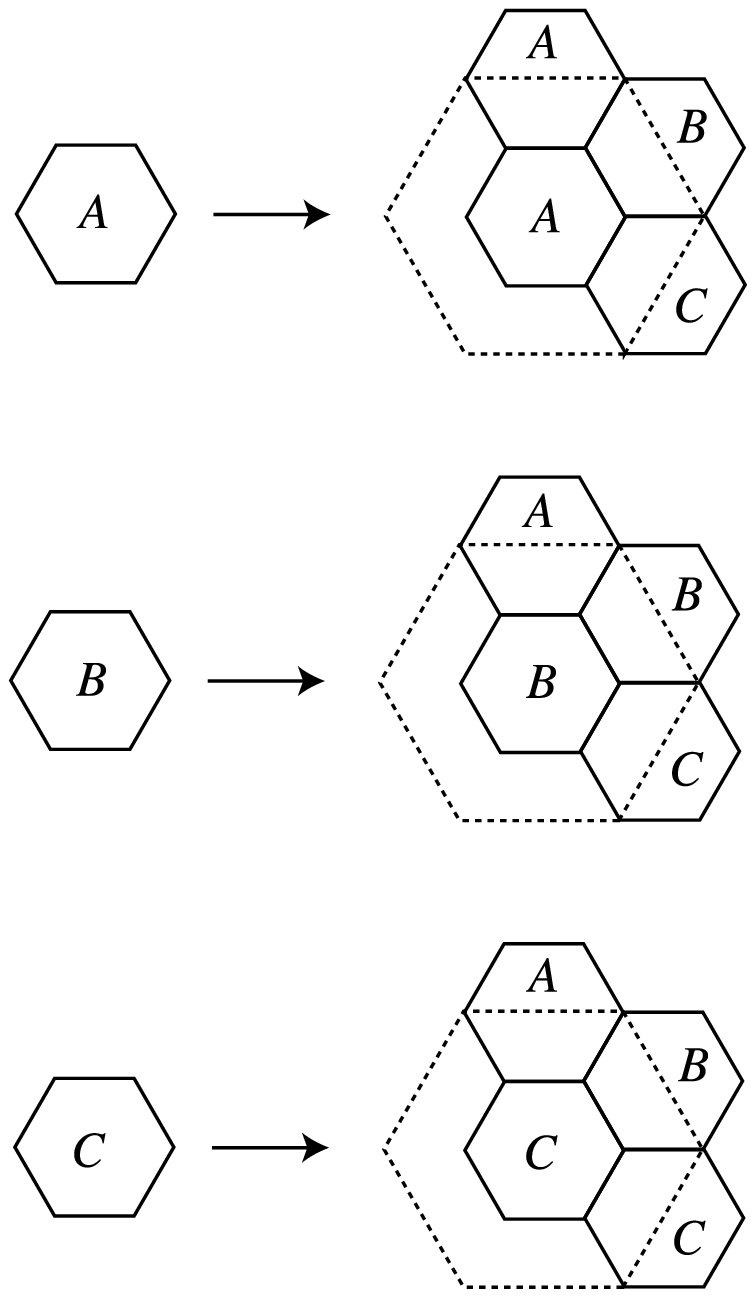}}
\caption{A substitution with amalgamation}
\end{figure}

In this section we assume that $\calT$ is a nonperiodic substitution tiling
space derived from a substitution-with-amalgamation $\sigma$. Since $\sigma$ 
maps $\calT$ to $\calT$, it induces a pullback map: $\sigma^*: H^1(\calT, \R^d)
\mapsto  H^1(\calT, \R^d)$. We decompose $H^1(\calT, \R^d)$ into (generalized)
eigenspaces of $\sigma^*$. The largest eigenvalue, denoted $\lambda_{PF},$ 
equals the stretching factor and has multiplicity $d^2$. (In terms of 
P-equivariant cohomology, a basis for this space is $dz^i \otimes e_j$,
where $e_1, \ldots, e_d$ is the standard basis for $\R^d$ and 
$z^1, \ldots z^d$ are Cartesian coordinates.) The corresponding eigenspace
is denoted $PF(\calT)$.  
The span of the (generalized) eigenvectors with eigenvalues
of magnitude strictly less than 1 is denoted $S(\calT)$.

Anderson and Putnam's results [AP] were originally stated for
substitutions without amalgamation, but apply also to
substitutions with amalgamation.  The substitution $\sigma$ induces a map
from $\Gamma^{(1)}$ to itself, and $\calT$ is the inverse limit of the
sequence of maps:
\begin{equation}
\Gamma^{(1)} \overset{\sigma}{\leftarrow} \Gamma^{(1)} 
\overset{\sigma}{\leftarrow} \Gamma^{(1)} 
\overset{\sigma}{\leftarrow} \Gamma^{(1)} \cdots
\end{equation}
The cohomology $H^1(\calT, \R^d)$ is then the direct limit of
$H^1(\Gamma^{(1)}, \R^d)$ under $\sigma^*: H^1(\Gamma^{(1)}, \R^d) \mapsto
H^1(\Gamma^{(1)}, \R^d)$.  Nonzero eigenspaces of $\sigma^*$ applied to
$H^1(\calT, \R^d)$ correspond to nonzero eigenspaces of $\sigma^*$ applied
to $H^1(\Gamma^{(1)}, \R^d)$.  The only difference between $H^1(\calT, \R^d)$
and $H^1(\Gamma^{(1)}, \R^d)$ is that $H^1(\Gamma^{(1)}, \R^d)$ may contain
a (generalized) zero-eigenspace of $\sigma^*$, while $H^1(\calT, \R^d)$
does not. 

This decomposition into eigenspaces of $\sigma^*$ makes it easy to 
identify the asymptotically negligible classes.

\begin{thm}\label{S=AN}
If $\calT$ is a substitution tiling space, then $N(\calT) = S(\calT)$.
\end{thm}

\demo{Proof.} For 
simplicity, suppose that the action of $\sigma^*$ on $H^1(\calT, \R^d)$
is diagonalizable.  
First we show that $N(\calT) \subset S(\calT)$.
For $\beta \in N(\calT)$, decompose $\beta$ as
\begin{equation}\label{breakbeta}
\beta = \beta_1 + \beta_2 + \cdots + \beta_k,
\end{equation}
with all terms nonzero, and with each $\beta_j$ an eigenvector of $\sigma^*$
with eigenvalue $\lambda_j$, with $|\lambda_1| \ge |\lambda_2| \ge \cdots
\ge |\lambda_k|$. Since $\beta_1 \ne 0$, by Theorem \ref{span-H_1}
there exists a recurrence class
$[p_z]$ such that $\beta_1[p_z] \ne 0$. The size of the recurrence 
$\sigma^n(p_z)$ grows exponentially with $n$, so $\beta(\sigma^n(p_z))$ 
goes to zero as $n \to \infty$. However,
\begin{equation}
\beta(\sigma^n(p_z)) = ((\sigma^*)^n \beta) p_z =
\lambda_1^n \beta_1(p_z) + \lambda_2^n \beta_2(p_z) + \cdots + \lambda_k^n
\beta_k(p_z).
\end{equation}
Since $\beta_1(p_z) \ne 0$, the only way this can converge to zero is
if $|\lambda_1|<1$, hence $\beta \in S(\calT)$. 

To show that $S(\calT) \subset N(\calT)$, note that every recurrence class
$[p_z]$ can be written as a sum
\begin{equation}
[p_z] = \sum_{j=0}^\infty (\sigma_*)^j[p_j],
\end{equation} 
with only a finite number of nonzero terms, and with all the $[p_j]$'s
belonging to a fixed bounded subset $U$ of $H_1(\Gamma^{(1)})$.  (The
argument is essentially that found in \cite[Proof of Lemma 3.2]{So2} 
and is not repeated here.) 
If the recurrence has size greater than 
$D(1+\lambda_{PF}+\cdots+\lambda_{PF}^{n-1})$
where $D$ is the recognition length, then the first $n$ terms are
zero.

Suppose that $\beta$ is an eigenvector of $\sigma^*$ with eigenvalue
$\lambda$ of magnitude strictly less than 1.  Let $K$ be larger than
the greatest value that $\beta$ takes on $U$. Then for any recurrence class
$[p_z]$, $|\beta(p_z)| < K/(1-|\lambda|)$.  If $p_z$ has size greater than
$D(1+\lambda_{PF}+\cdots+\lambda_{PF}^{n-1})$, then $|\beta(p_z)| < 
K |\lambda|^n/(1-|\lambda|)$.
Thus $\beta$ is asymptotically negligible.  Since linear combinations of 
asymptotically negligible classes are asymptotically negligible, $S(\calT)
\subset N(\calT)$. 

When $\sigma^*$ is not diagonalizable, the analysis is only slightly
more complicated. The expansion of $\beta(\sigma^n(p_j))$ may involve
$n$-th powers of eigenvalues times polynomials in $n$, rather than
just $n$-th powers of eigenvalues.  However, the conclusions are
unchanged.  \QED

\begin{thm}
If, for some non-negative integer $k$, either 
$\S(g) - (\sigma^*)^k \S(f) \in S(\calT)$ or $\S(f)-(\sigma^*)^k\S(g)
\in S(\calT)$, then $\calT_f$ is topologically conjugate to $\calT_g$.
\label{sufficient0}
\end{thm}

\demo{Proof.} The case where $k=0$ was already proven by Theorems
\ref{negligible} and \ref{S=AN}. All that remains is to prove
conjugacy when $\S(g) = \sigma^*(\S(f))$. Since the MLD class of
$\calT_g$ depends only on $\S(g)$, we can assume that $g = \sigma^*(f)$ as
cochains on $\Gamma^{(1)}$. 
[Note: since we need to work with $\Gamma^{(1)}$
rather than $\Gamma^{(0)}$, we pull our shape parameters $f$ and $g$
from $\Gamma$ back to $\Gamma^{(1)}$, and henceforth view them as
closed cochains in $C^1(\Gamma^{(1)}, \R^d)$.]   
That is, a cluster of tiles $\sigma(t_i)$ in
$\calT_g$ has the same dispacements as the single tile $t_i$ in $\calT_f$. 
To convert a tiling $x \in \calT_f$ to a tiling $\phi(x) \in \calT_g$,
replace each tile $t \in x$ by the cluster $\sigma(t)$, with shapes given
by $g$.  \QED

Theorem \ref{sufficient0} gives sufficient conditions for two tiling
spaces to be conjugate. For infinitesimal deformations, these conditions
are also necessary:

\begin{thm}[Theorem \ref{converse}]
For each shape parameter $f$ there is a neighborhood $U_f$ of $\S(f)$
such that, if $\S(g) \in U_f$ and $\calT_f$ and $\calT_g$ are 
topologically conjugate,
then $\S(g)-\S(f) \in S(\calT)$. 
\label{necessary0}
\end{thm}

The remainder of this section is devoted to proving Theorem \ref{necessary0}.
To prove the theorem we must
understand the extent to which patterns in $\calT$ 
can repeat themselves.  For each recurrence $(\z_1,\z_2)$, let the
\textit{degree} of the recurrence be its size divided by $|\z_2-\z_1|$.
For fixed $p>0$, we
show that the recurrence classes corresponding to recurrences with degree
$p$ or greater can be grouped into a
finite number of families, and we construct a conjugacy invariant from
the asymptotic displacements of these recurrences.  
Comparing these invariants for
different shape parameters then gives necessary conditions for conjugacy.

It is convenient to work with explicit matrices.  The set of integer
linear combinations of recurrence classes of size at least $2A$ 
(where $A$ is the diameter of the largest tile) is a sub-lattice of 
$H_1(\Gamma^{(1)})$. Let $\{a_1, \ldots, a_s\}$ be a basis for
this sub-lattice.  Since $\sigma$ maps recurrences to
recurrences, it maps the lattice to itself.  Let $M$ be the matrix
of this map relative to the basis $\{a_1, \ldots, a_s\}$.  Note that the
entries of $M$ are integers.  The matrix $M$ will play a role similar to
that of the substitution matrix in one-dimensional substitutions.  For
each recurrence class in $H_1(\Gamma^{(1)})$, 
the corresponding \textit{recurrence vector} in $\Z^s$ is the decomposition 
of the class in the $\{a_1, \ldots, a_s\}$ basis. 
For any shape parameter $f$, let $\L_f=(f(a_{1}),\ldots ,f(a_{s}))$ be
the ($\R^d$-valued row) vector that gives the displacements, in $\R^d$,
corresponding to lifts of the various loops.  For a recurrence with vector
$\v$, the corresponding displacement is $\L_f \v$.  The
Euclidean length of $\L_f \v$ is called the {\em length} of $\v$, and
denoted $|\v|_f$. We call $\L_f$ the
\textit{shape vector} of $\calT_f$.

Let $\calT_0$ be a tiling space whose shape vector $\L_0$ is a left
Perron-Frobenius eigenvector of $M$.  This can be the original substitution
tiling space, a linear transformation applied to the original space, 
or a space that is MLD to such a linear transformation.  In all such cases,
there exists a pseudo-self-similar tiling in the space, and we call
the space itself pseudo-self-similar. 

If $\v$ is a recurrence vector of degree $p$ in $\calT_0$,
corresponding to a recurrence $(\z_1,\z_2)$ in a tiling $x$, then $M\v$
is also a recurrence vector of degree at least $p$, corresponding to
the recurrence $(\lambda_{PF}\z_1, \lambda_{PF}\z_2)$ in the tiling
$\sigma(x)$. Both the size of the matching balls around $\z_1$ and
$\z_2$ and the distance from $\z_1$ to $\z_2$ get stretched by the same
factor $\lambda_{PF}$ in the substitution. Thus each recurrence vector
of degree $p$ gives rise to a family of recurrence vectors
$M^{k}\v$. The following theorem limits the number of such families.

\begin{thm}
Let $\calT_0$ be any pseudo-self-similar tiling space and let $p>0$. 
There is a finite collection of vectors $\left\{ \v_{1},\ldots
,\v_{N}\right\} $ such that every recurrence vector of degree at least
$p$ for $\calT_{0}$ is of the form $M^{k}\v_{i}$ for some pair $(k,i)$. 
\label{families}
\end{thm}

\medskip \noindent \textbf{Proof.\enspace}
As before, let $D$ be the recognition length of the substitution $\sigma$
acting on $\calT_0$. 
Suppose $\v$ is a recurrence vector of degree at least 
$p$ for $\mathcal{T}_{0}$,
where the radius $R$ of the matching balls is much greater than $D$. 
Then the supertiles within balls of radius $R-D$ around $\z_1$ and $\z_2$
must also agree.  Thus there is a recurrence vector 
$\v_{1}$ such that $\v=M\v_{1}$, and such that $\v_{1}$ is a recurrence vector
of degree at least $(R-D)/(\lambda_{PF}|\v_1|_0) \ge 
(p|\v|_0-D)/(\lambda _{PF}|\v_1|_0)
= p - (D/|\v|_0)$. 

Repeating the process, we find recurrence vectors $\v_i$ such that 
$M\v_i = \v_{i-1}$ and $\v_i$ is a recurrence vector of degree at least
\begin{equation}
p - \frac{D}{|\v_i|_0} (\lambda_{PF}^{-1} + \lambda_{PF}^{-2} + \cdots +
\lambda_{PF}^{-i}) \ge p - \frac{D}{|\v_i|_0} \sum_{\ell=1}^\infty
\lambda_{PF}^{-\ell} = p - \frac{D}{|\v_i|_0 (\lambda_{PF}-1)}.
\end{equation}

Now pick $\epsilon <p$. We have shown that every recurrence vector of
degree at least 
$p$ is of the form $M^{k}\v_{i}$, where $\v_{i}$ is a recurrence vector
of degree at least $p-\epsilon $, and where $|\v_{i}|_0$ is bounded by 
$\frac{
D\lambda _{PF}}{\epsilon (\lambda _{PF}-1)}$. However, 
by finite local complexity, there are only a
finite number of recurrence vectors of this length.\QED

\medskip

As defined, the degree of a recurrence vector depends on the shape vector
of the tiling space.  Indeed, since linear transformations do not preserve
lengths or ratios of lengths, the degree of a
recurrence vector in two pseudo-self-similar tiling spaces may not be the same.
However, for any tiling space $\calT_f$, the degrees of large recurrence
vectors are approximately the same as that of some pseudo-self-similar
tiling space: 

\begin{lem}
Let $\L_f$ be decomposed as $\L_f= \L_0 + \L_r$, where $\L_r$ is a linear
combination of (generalized) eigenvectors with eigenvalue less than
$\lambda_{PF}$ in magnitude, and $\L_0$ is a Perron-Frobenius eigenvector. Let 
$\calT_0$ be a pseudo-self-similar tiling space with shape vector $\L_0$. 
For each recurrence vector $\v$, let $r_\v = |\v|_f/|\v|_0$. 
For each $\epsilon >0$ there exists a length $R$ such that all 
recurrence vectors $\v$ with
$|\v|_f>R$ have $|r_\v - 1| < \epsilon$. \label{lem3.3}
\end{lem}

\medskip 

If $\lambda_2$ is the second largest eigenvalue in magnitude, then 
in recurrences of the form $M^n \v$, $|r_\v - 1|$ goes to 0 as 
$\left( \lambda _{2}/\lambda _{PF}\right) ^{n}$. Moreover, this convergence
is uniform in $n$ for all recurrences, because all non-zero recurrences  are longer than the diameter of a
ball that can fit in each prototile.
As indicated in the proof of theorem \ref{S=AN}, given any 
$N$,  all recurrences of 
sufficient length are sums of recurrences of the form $M^n \v$ with 
$n>N$, and so the result follows.
\QED
\medskip

\begin{lem}
Let $\epsilon >0$, and let $f$ be fixed. There is a self-similar tiling
space $\calT_0$ and a 
length $ R$ such that every recurrence vector of degree at least $p$ in
$\calT_f$, of length greater than $R$, 
is a recurrence vector of degree at least $p-\epsilon$ in
$\calT_0$. Furthermore, every recurrence vector of degree at least
$p+\epsilon$ in $\calT_0$, of length greater than $R$, is a recurrence 
vector of degree at least $p$ in $\calT_f$. 
\label{samerepvecs}
\end{lem}

\demo{Proof.} Decompose $\L_f=\L_0 + \L_r$ as above, and let $\calT_0$
be the pseudo-self-similar tiling space with length vector $\L_0$. By Lemma
\ref{lem3.3}, the ratio of lengths in $\calT_f$ and $\calT_0$
approaches 1 for large recurrences.  Thus the ratio of the radii of
the balls around $\z_1$ and $\z_2$ and the distance $|\z_2-\z_1|$ are
within $\epsilon$ for large recurrences.  \QED

Theorem \ref{families} showed there there are only a finite number of
families of recurrence vectors of a given degree for a self-similar
tiling.  Lemma \ref{samerepvecs} extends that result to all tiling spaces.  

We now construct a topological invariant from the asymptotic displacements
$\L_f \v$ of large recurrence vectors $\v$. 

\begin{thm}\label{staysame}
Suppose that $\phi: \calT_f \mapsto \calT_g$ is a topological conjugacy. Given
any positive constants $p, \epsilon_1, \epsilon_2$, there exists a positive
constant $R$ such that, for each recurrence vector $\v$ of degree at least
$p$ for $\calT_g$ with $|\v|_f>R$, there exists a recurrence vector 
$\v'$ of degree at least $p-\epsilon_1$ for $\calT_g$ with 
$|\L_f \v-\L_g \v'|<\epsilon_2$.
\label{invariant}
\end{thm}

In other words, up to small errors in the degree and the displacement that
vanish in the limit of large recurrence classes, the degrees and displacements
of the recurrence classes are conjugacy invariants. 

\demo{Proof.} Let $A$ be the diameter of the largest tile in the $\calT_g$ 
system. Since $\phi$ is uniformly continuous (being continuous with a 
compact domain), there exists a constant
$D_0$ such that, if $x$ and $y$ are tilings in $\calT_f$ that agree on a 
ball of radius $D_0$ around the origin, then $\phi(x)$ and $\phi(y)$ agree
on a ball or radius $A$ around the origin, up to a small translation.  Thus
if $x$ and $y$ agree on a ball of radius $R>D_0$, then $\phi(x)$ and $\phi(y)$
agree on a ball of radius $R-D_0$, up to a translation whose norm is bounded
by a decreasing function $h(R)$, with $\lim_{R \to \infty} h(R)=0$. 

Now let $\v$ be a recurrence vector of degree $p$ for $\calT_f$, representing
a recurrence $(\z_1,\z_2)$ in a tiling $x \in \calT_f$.
Then $x-\z_1$ and $x-\z_2$ agree on a ball of radius $p|\v|_f$ about
the origin, so $\phi(x)-\z_1$ and $\phi(x)-\z_2$ agree on a ball of
radius $p|\v|_f-D_0$, up to translation of size at most $h(R)$. 
Thus there exists a point $\z_3$, within $h(D)$ of $\z_2$, such that 
$\v'=[(\z_1,\z_3)]$ is a
recurrence class of degree at least 
$\frac{p|\v|_f-D_0}{ |\v|_f+h(D)}$ in $\calT_g$.
For $D$ large enough, this is greater than $p-\epsilon_1$ and 
$|\L_f \v- \L_g \v'| < h(R)$ is less than $\epsilon_2$. 
\QED

\demo{Proof of Theorem \ref{necessary0}.}

Note that all recurrence vectors are recurrence vectors of some positive
degree.  We can therefore pick $p_0$ such that the integer span of 
the recurrence vectors of degree at least 
$p_0$ is an $s$-dimensional sublattice of $\Z^s$.  Pick $\epsilon$
small enough (and adjust $p_0$ by up to $\epsilon$, if necessary) so 
that all the families of recurrence vectors of degree at least
$p_0-\epsilon$ are
also families of recurrence vectors of degree at least 
$p_0+\epsilon$. Let $\v_1,
\ldots, \v_k$ be generating vectors of those families. By multiplying by
appropriate powers of $M$, these can be chosen so that all of the magnitudes
of the displacements $|M^k \v_i|_f$ are within a factor of $\lambda_{PF}$
of one another for large $k$. 

Now pick a neighborhood $U_\epsilon$ of $f$ in $C^1(\Gamma, \R^d)$
small enough that, if $g' \in U_\epsilon$, then every recurrence class 
in $\calT_f$ of degree at least
$p_0+\epsilon$ is also a recurrence class in $\calT_{g'}$ of
degree at least $p_0$, and such that every recurrence class in 
$\calT_{g'}$ of degree at least $p_0$
is a recurrence class in $\calT_f$ of degree at least $p_0-\epsilon$. 
This insures that
the families of recurrence classes of degree at least 
$p_0$ in the two tiling spaces are exactly 
the same. 

If $\calT_f$ and $\calT_{g'}$ are conjugate then, by Theorem
\ref{invariant}, the displacements $ \L_f M^k \v_j$ can be approximated by
$\L_{g'} M^{k'} \v_{j'}$ for some $k', j'$, and the approximation get
successively better as $k \to \infty$.  However, if $U_\epsilon$ is
chosen small enough, the only values of $k',j'$ that come close to
approximating are $k'=k$ and $j'=j$.  Thus the limit of $(\L_f-\L_{g'})
M^k \v_j$
must be zero.  By Theorem 2.3, this implies that $(\sigma^*)^k(\S(f)-\S(g'))$
approaches zero as $k \to \infty$, and hence that 
$\S(f)-\S(g') \in S(\calT)$.

Finally, let $U_f = \S(U_\epsilon)$.
If $\S(g) \in U_f$, then $\calT_g$ is MLD to a tiling space $\calT_{g'}$
with $g' \in U_\epsilon$ and $\S(g')=\S(g)$. Since $\calT_{g'}$ is
conjugate to $\calT_f$, $\S(f)-\S(g) = \S(f)-\S(g') \in S(\calT)$.
\QED




\section{Ergodic Properties}

Now we turn to the topological point spectrum of substitution tiling
spaces, by which we mean the eigenvalues of continuous eigenfuctions
of the translation action.  First we determine general criteria for a
vector to be an eigenvalue for a continuous eigenfunction, and then we
apply this criteria to some special cases, depending on the form of
the matrix $M$ defined in Section \ref{sub}. This eventually leads to
criteria for topological weak mixing.
 
\begin{thm}
The vector $\k\in\R^d$ is in the point spectrum of
$\mathcal{T}_{f}$ if and only if, for every recurrence vector
$\v$,
\begin{equation}\label{eig}
\dfrac{1}{2\pi } (\k \cdot \L_f) M^{m}\v\rightarrow 0 \text{(mod 1) as } 
m\rightarrow
\infty , 
\end{equation}
where the convergence is uniform in the size of $\v$. 
%

\label{zeromod1}\ \end{thm} \medskip

\demo{Proof.}  Let $x_0$ be a tiling in $\calT_0$ fixed by the substitution
homeomorphism, let $x$ be its image in $\calT_f$ under the homeomorphism
of Theorem 1.1, and let $E:\mathcal{T}_{f} \mapsto S^1$ be a continuous
eigenfunction with eigenvalue $\k$.  Let $\v$ be a recurrence vector,
then there is a recurrence $(\z_1,\z_2)$ in $x$ with vector $\v$ of
some size $s_0$. By applying the substitution homeomorphism $m$ times,
we obtain a recurrence $(\z_1^m,\z_2^m)$ 
in $x$ with recurrence
vector $M^{m}\v$, whose displacement is $\z_2^m - \z_1^m = \L_f M^m \v$.  
Then $x-\z_1^{m}$ and $x-\z_2^{m}$ agree on
patches of size $s_m$, where $s_m \rightarrow \infty$ as $m
\rightarrow \infty.$ Hence
\begin{equation}
1=\text{lim}_{m \rightarrow \infty} \dfrac{E(x-\z_2^{m})}{E(x-\z_1^{m})}
=\text{lim}_{m \rightarrow \infty} \dfrac{E(x)\text{exp}(-i\k \cdot \z_2^{m})}
{E(x)\text{exp}(-i\k \cdot \z_1^{m})}=
\text{lim}_{m \rightarrow \infty} \text{exp}(-i\k\cdot (\z_2^{m}-\z_1^{m})).
\end{equation}
Thus, we obtain Equation \ref{eig}, and the uniform convergence follows
from the uniform continuity of $E$.

Conversely, assume that we have the stated convergence for all
recurrence vectors $\v$ for some $\k \in \R^{d}$. We construct a
continuous eigenfunction $E$ by first assigning $x$ the value $1$.
Then we necessarily have for any $\z \in \R^{d}$,
$E(x-\z)=\text{exp}(-i\k \cdot \z)$. To show that $E$ extends as required to
all of $\mathcal{T}_{f}$, it suffices to show that $E$ as so defined
is uniformly continuous on the orbit of $x$. But, given an $\epsilon
>0$, by Equation \ref{eig} if $x-\z_1$ and $x-\z_2$ agree on patches
of sufficiently large size up to a small translation, we have that
$E(x-\z_1)$ and $E(x-\z_2)$ agree to within $\epsilon$.
\QED

The application of this criterion depends on the eigenvalues of $M$
and on the possible forms of the recurrence vectors.

\begin{thm}
Suppose that all the eigenvalues of $M$ are of magnitude 1 or greater.
If $\k$ is in the spectrum, then all elements of $\k \cdot \L_f/2\pi$ are
rational.  
\label{allbig1}
\end{thm}

\medskip \noindent \textbf{Proof.\enspace} Let $\k$ be in the point
spectrum, and consider the sequence of real numbers
$t_{m}=(\k \cdot \L_f)M^{m}\v/(2\pi )$, where $\v$ is a fixed recurrence
vector. Let $p(\lambda )=\lambda ^{n}+a_{n-1}\lambda ^{n-1}+\cdots
+a_{0}$ be the characteristic polynomial of $M$. Note that the
$a_{i}$'s are all integers, since $M$ is an integer matrix. Since
$p(M)=0$, the $t_{m}$'s satisfy the recursion:
\begin{equation}
t_{m+n}=-\sum_{k=0}^{n-1}a_{k}t_{m+k}.  \label{intrecur}
\end{equation}
By Theorem \ref{zeromod1}, the $t_{m}$'s converge to zero $\pmod{1}$. 
That is, we can write 
\begin{equation}
t_{m}=i_{m}+r_{m}  \label{splits}
\end{equation}
where the $i_{m}$'s are integers, and the $r_{m}$'s converge
to zero as real numbers. By substituting the division (\ref{splits})
into the recursion (\ref{intrecur}), we see that both the $i$'s and
the $r$'s must separately satisfy the recursion
(\ref{intrecur}), once $m$ is sufficiently large. However, any
solution to this recursion relation is a linear combination of powers
of the eigenvalues of $M$ (or polynomials in $m$ times eigenvalues to
the $m$-th power, if $M$ is not diagonalizable).  Since the
eigenvalues are all of magnitude one or greater, such a linear
combination converges to zero only if it is identically
zero. Therefore $r_{m}$ must be identically zero for all
sufficiently large values of $m$.

Apply this procedure to $s$ linearly independent recurrence vectors
$\v_1,\ldots,\v_s$, and pick $m$ large enough that the corresponding
$t_m(\v_i)$ are integers for each $i=1,\ldots,s.$ Note that $t_{m}$ is
an integer linear combination of the elements of the vector $ \k \cdot
\L_f/(2\pi )$. However,
\begin{equation}
(t_m(\v_1),\ldots,t_m(\v_s)) = (\k \cdot \L_f/2\pi) M^m (\v_1,\ldots, \v_s)
\end{equation}
The matrices $M$ and $V= (\v_1,\ldots,\v_s)$ are invertible and have integer 
entries, so by Cramer's rule their inverses have rational entries. 
Thus the components of $\k \cdot \L_f/(2\pi )$ must all be rational. 
\hfill $\square$\medskip

\begin{cor}
Suppose all the eigenvalues of $M$ have magnitude 1 or greater. Let $G
\subset \R^d$ be the free Abelian group generated by the entries of $\L_f$.
Let $\bar G$ be the closure of $G$ in $\R^d$, and let $P$ be the identity
component of $\bar G$. The point spectrum of $\calT_f$ lies in the orthogonal
complement of $P$. 
\end{cor}

In dimension greater than 1, the assumptions of Theorem \ref{allbig1} 
are rarely met. When eigenvalues of magnitude less than 1 exist, the 
conclusions are somewhat weaker.

\begin{thm}
\label{somesmall1} 
Let $S$ be the span of the (generalized left-) eigenspaces of $M$ with
eigenvalues of magnitude strictly less than 1. If $\k$ is in the point
spectrum, then $\k \cdot \L_f/ {2 \pi}$ is the sum of a rational vector and an
element of $S$.
\end{thm}

\medskip\noindent\textbf{Proof.\enspace} 
Let $V$ be the matrix $(\v_1,\ldots,\v_s)$ as in the proof of 
Theorem \ref{allbig1}.  Construct the row vector 
\begin{equation}
\t_m = \dfrac{1}{2\pi } (\k \cdot  \L_f) M^m V= 
\dfrac{1}{2\pi } (\k \cdot \L_f) V (V^{-1}MV)^m.
\end{equation}

As in the proof of Theorem \ref{allbig1}, each entry of $\t_m$
converges to zero (mod 1), so we can write $\t_m = \i_m +  
\r_m$,
with each entry of $\i_m$ integral and 
$\r_m$ converging to zero, and with the eventual conditions
\begin{equation}
\i_{m+1} = \i_m V^{-1} M V; \hspace{.1in}
\r_{m+1} = \r_m V^{-1} M V.
\end{equation}
Since the $\r_m$'s converge to zero,
they must lie in the span of the 
small eigenvalues of $V^{-1} M V$.  
Thus, by adding an element of $S$ to $(\k \cdot  \L_f)/2\pi$, 
we can then get all the $\r_m$'s to be identically
zero beyond a certain point.  Since $M$ is invertible, 
the resulting value of $\k \cdot \L_f/2\pi$ must then be
rational. \QED
\medskip

Another way of stating the same result is to say that $(\k \cdot \L_f)/2\pi$,
projected onto the span of the large eigenvectors, equals the
projection of a rational vector onto this span.

This theorem can be used in two different ways. First, it constrains
the set of shape parameters (that is, vectors $\L_f$) for which the system
admits point spectrum. Let $d_b$ be the number of large eigenvalues,
counted with (algebraic) multiplicity.  There are only a countable
number of possible values for the projection of $\k \cdot \L_f/2\pi$ onto
the span of the large (generalized) eigenvectors. In other words, one
must tune $d_b$ parameters to a countable number of possible values
in order to achieve a point in the spectrum.  Of course, $\k$ itself gives
$d$ parameters.  Thus we must tune at least $d_b-d$ additional parameters to
have any point spectrum at all.  In particular, if $d_b>d$, then a generic
choice of shape parameter gives topological weak mixing, proving Theorem
\ref{mixing}. 
(Note that the Perron-Frobenius
eigenvector always occurs with multiplicity $d$, so that $d_b-d$ is never
negative.) 

A second usage is to constrain the spectrum for fixed $\L_f$.  The
rational points in $\red^s$, projected onto the span of the large
eigenvalues, and then intersected with the $d$-plane defined by a fixed $\L_f$
(i.e., the set of all possible products $(\k \cdot \L_f)/2\pi$),
forms a vector space over $\qed$ of dimension at most $s+d-d_b$. As a
result, the point spectrum tensored with $\mathbb{Q}$ is a vector
space over $\qed$ whose dimension is bounded by $d$ plus the number of
small eigenvalues.  Below we derive an even stronger result, in which
only the small eigenvalues that are conjugate to the Perron-Frobenius
eigenvalue contribute to the complexity of the spectrum.

\begin{thm}
\label{notfull} Let $b_{PF}$ be the number of large eigenvalues, 
counted without multiplicity, that
are algebraically conjugate to the Perron-Frobenius eigenvalue
$\lambda_{PF}$ (including
$\lambda_{PF}$ itself), and let
$s_{PF}$ be the number of small eigenvalues conjugate to $\lambda_{PF}$. 
For fixed $\L_f$, the dimension over $\qed$ of the point
spectrum tensored with $\qed$ is at most $d(s_{PF}+1)$.  
\end{thm}

\medskip\noindent\textbf{Proof.\enspace} As a first step we
diagonalize $M$ over the rationals as far as possible.  By rational
operations we can always put $M$ in block-diagonal form, where the
characteristic polynomial of each block is a power of an irreducible
polynomial.  Since the Perron-Frobenius eigenvalue $\lambda _{PF}$ has both
geometric and algebraic
multiplicity $d$, every eigenvalue algebraically conjugate to $\lambda
_{PF}$ also has multiplicity $d$. Thus there are $d$ blocks whose
characteristic polynomial has $\lambda _{PF}$ for a root.   We
consider the constraints on the spectrum that can be obtained
from these blocks alone.

Consider the projection of $\k \cdot \L_f/(2\pi)$ onto the large
eigenspaces of the Perron-Frobenius block.  Since only $d(b_{PF}+s_{PF})$
components of $\t$ (expressed in the new basis) contribute, this
is the projection of $\qed^{d(b_{PF}+s_{PF})}$ onto $\red^{d b_{PF}}$,
whose real span is all of $\red^{d b_{PF}}$.  Intersected with the $d$-plane
defined by a fixed $\L_f$, this gives a vector space of dimension at most
$d(s_{PF}+1)$ in which $\k$ can live.  \hfill $\square$

\section{One dimensional substitutions revisited}

When discussing 1-dimensional substitutions, with $n$ tile types $t_1,
\ldots, t_n$, the conventional object of study is the {\em
substitution matrix}, whose $(i,j)$ entry gives the number of times
that $t_i$ appears in $\sigma(t_j)$.  Indeed, in our previous study
\cite{CS} of one dimensional tilings, all the results were phrased in
terms of eigenvalues and eigenspaces of the substitution matrix,
rather than on the induced action of $\sigma$ on homology.  These
results become much simpler when viewed homologically.  In this section,
the substitution matrix will be denoted $M_s$, while the matrix that
gives the action of $\sigma$ on a basis of recurrences will be denoted
$M_h$.

In \cite{CS} we defined, for each recurrence $(z_1,z_2)$, a vector in
$\zed^n$ that listed how many of each tile type appears in the
(unique) path from $z_1$ to $z_2$. This vector $\v$ was called {\em
full} if the vectors $(M_s)^k \v$, with $k$ ranging from 0 to $n-1$,
were linearly independent.  Many of our theorems required the
existence of a recurrence with a full vector. This is a strong
condition, as it implies that $H_1(\Gamma)$ is a lattice of rank $n$.
This is true when $n=2$, or when the characteristic polynomial of
$M_s$ is irreducible, but is typically false for more complicated
substitutions.  When $H_1(\Gamma)$ has rank less than $n$, there are
deformations of tile lengths that have no effect on the lengths of
recurrences, and so lead to MLD tilings.  By looking at $M_h$ rather
than $M_s$, we automatically avoid those extraneous modes.

Consider the difference between the following theorem, proved in \cite{CS},
and its restatement in terms of $M_h$:

\begin{thm}[CS] Suppose that all the eigenvalues of $M_s$ are of magnitude
1 or greater, and that there exists a recurrence with a full vector. If
the ratio of any two tile lengths is irrational, then the
point spectrum is trivial.
\end{thm} 

\begin{thm}[Corollary of Theorem 4.2] Suppose that all the 
eigenvalues of $M_h$ are of magnitude 1 or greater.  If the ratio of
the lengths of any two recurrences is irrational, then the point
spectrum is trivial.
\end{thm}

In addition to $M_h$ not containing irrelevant information found in $M_s$,
$M_h$ may contain some relevant information {\em not} found in $M_s$.
If $H_1(\Gamma^{(1)})$ has higher rank than $H_1(\Gamma^{(0)})$, then $M_h$
contains information about the dynamical impact of changing the sizes of
the collared tiles, and not merely the effect of changing the original,
uncollared tiles.  

As an example, consider the Thue-Morse substitution $(a \to ab, b \to
ba)$, in which $M_s=\begin{pmatrix} 1 & 1 \cr 1 & 1 \end{pmatrix}$ has
eigenvalues 2 and 0.  $\Gamma$ is the wedge of two circles, one
representing the tile $a$ and one representing the tile $b$, so
$H_1(\Gamma)=\zed^2$, and the action of $\sigma$ on $H_1(\Gamma)$ is
described by $M_s$.  However, $H_1(\Gamma^{(1)})$ has rank 3
\cite{AP}, and the eigenvalues of $M_h$ are 2, $-1$, and 0. The
additional large eigenvalue $-1$ shows that the dynamics of the
Thue-Morse tiling space are in fact sensitive to changes in tile size.
Changes in the size of the {\em uncollared} tiles have no qualitative
effect, but changes in the size of the collared tiles (i.e., changes
in tile size that depend on the local neighborhood of
those tiles) can eliminate the point spectrum.

\bigskip

\noindent{\em Acknowledgements.} We thank Michael Baake, 
Charles Holton, Johannes 
Kellendonk, Brian Martensen, Charles Radin and Bob Williams for helpful
discussions.

\end{document}